\documentclass[a4paper]{article}

\usepackage[english]{babel}
\usepackage[utf8]{inputenc}
\usepackage{amsmath}
\usepackage{graphicx}
\usepackage[colorinlistoftodos]{todonotes}

\title{The Efficacy of the Flipped Classroom Technique in Undergraduate Math Education: A Review of the Research}

\author{Adeli Hutton}
\date{}
\begin{document}
\maketitle
\begin{abstract}
     The flipped classroom technique has recently been a focus of attention for many math instructors and pedagogical researchers. Although research on the subject has greatly increased in recent years, it is still debated whether the flipped classroom technique can significantly increase the overall success of students in undergraduate math courses. While there have been meta-analyses that consider the effectiveness of the technique across university disciplines and within other STEM fields, there has not yet been a systematic review within undergraduate math education. By analyzing the existing research and compiling the quantitative and qualitative data, this paper examines the efficacy of the flipped classroom technique in undergraduate math courses, ranging from introductory calculus to transition-to-proof courses, in regards to students’ performance in the classes, perceptions of the technique, and associated self-efficacy. This paper also introduces the current use of the technique and covers successful implementation methods. Additionally, it highlights the flipped classroom technique’s potential for improving retention of members of underrepresented groups in math by increasing their sense of belonging as well as discusses the efficacy of the method in early proof-based courses in regards to students’ acquisition of sociomathematical norms.
\end{abstract}

\section{INTRODUCTION}

{\setlength{\parskip}{3mm}\subsection{Traditional Classroom Model}}
\hspace{4mm} The traditional pedagogical model of undergraduate math almost goes without explanation: professors prepare and present lectures to students, who listen, take notes, and ask questions. Student learning is typically assessed and graded through a combination of homework assignments, quizzes, and exams. This traditional lecture style remains the norm in undergraduate math education, especially in the U.S., and students are usually comfortable with this model of learning due to their familiarity with it [2]. Betty Love et al. who studied students’ learning and perceptions in a flipped classroom linear algebra course note that while the traditional method generates positive results in learning for many students, it is not an ideal learning approach for every student [2]. In an article studying learning and teaching styles in college science education, Richard Felder emphasizes that students can have significantly different learning styles. According to Felder, when students’ learning styles are compatible with the teaching style of an instructor, students tend to “retain information longer, apply it more effectively, and have more positive post-course attitudes toward the subject” than students who were in course settings that did not match their learning styles [3]. Love et al. note that the traditional lecture style is “well entrenched” in college instruction and is how most professors and faculty members were themselves taught. Though most professors use the traditional lecture method due to their familiarity with it, they do not necessarily think it is the best method for optimizing student learning [2]. A survey conducted by Susan Lord and Michelle Camacho, engineering faculty at the University of San Diego, found that only 36 percent of teaching-oriented engineering faculty think that the traditional lecture is a good teaching approach. However, 60 percent of the respondents still teach in this traditional method [4].
\subsection{Toward a New Model}  \hspace{4mm}  For a long time, college math faculty have noticed students’ disengagement and passivity during traditional lectures and have linked this to lower retention levels of students [5]. Felder mentions that many students who enter college with the initial intention and ability to pursue a degree in a science field switch to non-science areas of study after completing introductory science courses, due in part to the passive-learning nature of the classrooms of these courses. He hypothesizes that this subset of college students could be numerous enough to prevent a deficit of American scientists and engineers [3]. Research has shown that active learning helps to combat these challenges of student disengagement and passivity, which lead to disinterest in certain subject areas. In his book Flipped Learning: A Guide for Higher Education Faculty, Robert Talbert defines active learning as “any instructional method that engages students in the learning process, in an active way, as part of the group space activities” [6]. Research has shown that when professors give students opportunities to be more active during class, whether through participating in a group discussion or completing a worksheet, this yields increases in students’ learning outcomes [5]. In a meta-analysis of 225 studies reporting data on exam scores and failure rates in undergraduate STEM courses, Scott Freeman et al. found that students in traditional lecture courses had a 55 percent higher course failure rate than students in classes that included active learning and that exam grade averages in active-learning classes were half a letter grade (6 percent) higher than in traditional lecture classes. They note that these findings hold in analysis across STEM disciplines and with varying class sizes [7]. In an effort to increase active learning, some college math faculty have begun using a new instructional technique: the flipped classroom.

\subsection{The Flipped Class Model}
     \hspace{4mm} The idea of a flipped classroom, is that prior to and outside of class, students learn the material that would normally be covered in lecture, so time in class can be allocated to working through examples and problems. Cynthia Brame of Vanderbilt University’s Center for Teaching describes it as a model in which students acquire first-exposure learning before attending class, enabling them to focus on processing aspects of learning while in class [8]. Karjanto and Simon, who studied the effectiveness of the flipped classroom technique, note that there is no exact definition for flipped classrooms within research articles, however, for their study, they define the concept as doing “lecture at home” and “homework at school” [9]. The Flipped Learning Network defines flipped learning as “a pedagogical approach in which direct instruction moves from the group learning space to the individual learning space, and the resulting group space is transformed into a dynamic, interactive learning environment where the educator guides students as they apply concepts and engage creatively in the subject matter” [10]. The idea central to all definitions of the technique is that the method enables students to work through problems in class while there is an instructor and possibly a teaching assistant available to help them. During a traditional lecture, students typically do not have the opportunity to work problems or think about the concepts critically until they are outside of class, when there is not an instructor present for guidance and assistance. Love et al. mention that the current generation of students believes that “doing is more important than knowing and that learning is a trial-and-error process” and that many aspects of the flipped classroom, including the hands-on activities and increased interaction with instructors, may accommodate this learning style better than the traditional lecture method [2].
    {\setlength{\parskip}{3mm}
\section{IMPLEMENTATION OF THE METHOD}

\subsection{How the Flipped Model Works}

\hspace{4mm} Since there is no standardized definition for the concept of a flipped classroom, there are various ways and degrees to which the technique has been implemented. Instructors using the flipped classroom technique typically provide students with a pre-class lecture to watch at their convenience before going to class, assign reading from a chosen lesson in their textbooks, or both. The students are sometimes given a pre-class assignment or quiz over the material covered. The research studies emphasize the importance of providing an incentive for students to come to class prepared. In her article, Cynthia Brame outlines four key elements of the flipped classroom method [8]. Her key elements along with ideas from other articles are summarized: \\
\\
1) Students acquire first exposure learning prior to class.\\
\\
The medium for this varies and can include textbook reading and lecture videos. The course instructor can create video lectures or find them on Khan Academy, YouTube, MIT OpenCourseWare, or other online sources and assign them to the students. Other resources that may be provided to students prior to class are annotated notes, pre-readings, study guides, and automated tutoring systems [11]. Love et al. note that their study suggests that instructional videos are a key element in implementing flipped classrooms due to their popularity with students [2]. \\ \\
2) Students are provided with an incentive to prepare for class. }\\ \\
Students complete a task that is associated with points for their course grade. It can be an online quiz, worksheet, or any assignment that causes them to complete and think about the first exposure learning that they were assigned. \\ \\
3) Instructors create a way to assess students’ understanding. \\ \\
This can be accomplished through the pre-class assignment if instructors are able to view students’ results on it prior to class. This allows course instructors to pinpoint areas in which students are struggling and emphasize and clarify these concepts during class. This process is known as Just-in-Time Teaching and is a pedagogical strategy developed by Gregor Novak. Love et al. mentions that Just-in-Time Teaching is often implemented in flipped class settings [2]. \\ \\
4) Instructors structure in-class activities to focus on higher level cognitive activities.  \\ \\
\hspace{4mm} Because students learn the necessary basic knowledge outside of class, time in class can be devoted to activities that promote deeper learning. This is the key aspect of the flipped classroom technique—that students are able to use class time to increase their skill levels and deepen their understanding of mathematical content. In math courses that use the flipped classroom technique, professors typically choose or create problems for students to work through in class, often in groups. These problems are typically chosen to encourage critical thinking of concepts that professors want to emphasize. This gives students the opportunity to work problems that are often of a higher difficulty level while the teacher is available for questions. During class, ideally professors are also able to differentiate instruction to better suit the needs of individual students, as the skill level of students within a class can greatly vary. Robert Talbert mentions that in implementing the flipped classroom technique, he can formulate basic exercises to help groups of students who are struggling with the concepts as well as create new and related problems for groups of students who have already completed the worksheet for the day [12].\\ \\
    \hspace{4mm} The key elements of flipped learning that Brame mentions are very similar to the core features that O’Flaherty and Phillips discuss in the conclusion of their scoping review of the flipped classroom: “content in advance (generally the pre-recorded lecture), educator awareness of students’ understanding, and higher-order learning during class time” [11].
\subsection{Use of Technology in Flipped Classes}
\hspace{4mm} Many articles discuss the utilization of the advancement of technology as an important aspect of implementing the flipped classroom technique. Dominic McGrath et al. mention that technology can play an integral role in augmenting the advantages that the flipped classroom technique provides. It can create a way for instructors to interact with students outside of class as well as increase flexibility for watching recorded lectures at a chosen time and pace that best suits students’ learning of the content [13]. McGrath et al. also note that technology can allow for improved monitoring of students’ understanding and engagement within the course, helping to better identify students who are struggling [13]. Specifically, it can assist instructors in better applying the Just-in-Time Teaching technique if students are able to take quizzes or submit questions online to their instructors before class, allowing instructors to focus on areas that students need help with in class. For example, in their study on the effectiveness of the flipped classroom technique in an applied linear algebra course, Love et al. implemented daily readiness assessments into the flipped class section of the study [2]. These assessments were administered through Blackboard, were part of students’ grades, and consisted of three questions. Two of the questions related to the course content and the last question was always the following: “What did you find difficult or confusing about this section? If nothing was difficult or confusing, what did you find most interesting? Please be as specific as possible” [2]. Students were to submit their responses two hours before their class so that the instructor could read their responses prior to class to decide what to focus on during class [2].\\
     \\
\hspace{4mm} Technology can also increase the accessibility of communication between not only students and their instructors but also students and their peers. In his flipped classroom section of a transition-to-proof course, Robert Talbert created an online discussion forum using Piazza, a website that allows students to ask and answer questions of their course instructors and peers. He found that with Piazza, some students in his class became very active in helping their peers with the course material, and questions that were left unresolved on the forum became discussion topics for class [12]. Technology can also allow students to better view three-dimensional graphs, which is particularly helpful for courses in areas such as multivariable calculus and differential geometry, and instructors can incorporate these graphs into their video lectures or pre-class assignments. Although there are many benefits to using technology in any math class, McGrath et al. advise against assuming that the presence of increased technology in flipped classrooms is automatically productive. They note that a small percentage of students may struggle to use the technology with ease due to insufficient experience with it, and they emphasize the importance of providing comprehensive instructions and being sure that students have access to IT assistance [13].
\subsection{Development of the Flipped Classes}
 \hspace{4mm}  Williams et al. note that although the traditional lecture model has been “a staple of academia for close to a millennium,” the ideas behind the flipped classroom technique can be viewed as an even older pedagogical method in which class time focused more on discussion and academic debate than on direct transfer of information from the instructor to students [14]. Talbert emphasizes that individual instances of flipped learning could go back hundreds of years as professors may have structured their courses in this form [6]. Within published research, one of the earliest instances of flipped learning was the work of Eric Mazur, a Harvard University physics professor, in the 1990s. Mazur realized that his students were very skilled at successfully completing physics problems that involved computation and memorization but lacked conceptual understanding even after completing a semester of his physics class [6]. He remodeled his physics course and used team learning and in-class activities in place of lectures [14]. In order to ensure that students learned from their assigned reading, Mazur developed a software which allowed students to read the textbook online and ask questions and comment within the text so that other students and instructors could read and reply to what they wrote. Mazur’s course instruction design method has been known as peer instruction [6]. Talbert notes that Mazur’s example demonstrates that flipped learning “emerged as a solution to a concrete pedagogical problem associated with students’ conceptual understanding of a complex subject” [6]. Karjanto and Simon state that the term flipped classroom was “coined and popularized” by Jonathan Bergmann and Aaron Sams in 2007 [9]. They created video lectures to reverse, or flip, when lectures and homework took place and are believed to be the first to have done so [9]. Around 2000, Maureen Lage, Glenn Platt, and Michael Treglia, economics professors, discovered that students in introductory economics courses were struggling due to their learning styles not matching well with the traditional lecture method. In addition, professors teaching this course felt constrained by time in their ability to vary instructional methods [6]. Lage, Platt, and Treglia remodeled their classroom so that students would read and watch videos prior to class. The first ten minutes of class were allocated to addressing students’ questions over these pre-class assignments, and the remainder of the in-class time was spent on “lively activities,” review questions, and practice worksheets [6]. They named their format the inverted classroom. \\
 \\
\hspace{4mm} Many of these early examples of use of the flipped classroom technique were developed in part to improve students’ understanding of the more conceptual ideas of the classes. These ideas link to Bloom’s Taxonomy, a categorization system for educational goals developed by Benjamin Bloom, Max Englehart, Edward Furst, Walter Hill, and David Krathwohl in 1956 [15]. Since its development, Bloom’s Taxonomy has been applied to K-12 and college instruction alike and has evolved into the pyramid included below consisting of six categories: remember, understand, apply, analyze, evaluate, and create [15]. \\ \\
  {\setlength{\parskip}{3mm}
    \hspace{4mm}  In the context of the flipped classroom, students complete the lower tiers of cognitive work outside of class by gaining knowledge and comprehension through their pre-class assignments, which enables them to focus on the higher tiers of cognitive work—application, analysis, and knowledge synthesis—during class when their instructor is available [8]. In contrast, the traditional lecture method places emphasis on lower tiers of cognitive thinking during class time [16]. Karjanto and Simon note that in traditional lecture classrooms, students are on their own to complete these higher levels of thinking. In disciplines such as mathematics that have strong connections between course topics, this can cause students to be frustrated when they struggle in completing their homework on their own [9]. Karjanto and Simon refer to the framework of the flipped classroom as an inverted version of Bloom’s Taxonomy since higher order learning takes place in class and more time is allocated for it [9]. In terms of proof-based math classes, Robert Talbert mentions that the traditional lecture model has a practical flaw: students are expected to construct proofs—the most cognitively complex work—outside of the classroom. Because of the lack of immediate expert guidance in this model, students who struggle may give up or not be engaged in the class [12].
\subsection{Current Use of the Technique}
     \hspace{4mm} With increased interest in the method, the use of the flipped classroom technique is growing across higher education. Karabulut-Ilgu et al. mention that many universities have embraced flipped learning and cite a study that was conducted by the Center for Digital Education and Sonic Foundry that found that 29 percent of higher education faculty in the U.S. are currently implementing flipped learning and an additional 27 percent plan to implement it in the near future [17].
\subsection{Aspects of Successful Implementation}
     \hspace{4mm} Since there is no standardized model for the flipped classroom, instructors who decide to try the flipped classroom technique should be aware that proper implementation of the method is key for its success. Wasserman et al. mention that a potential concern with the growing trend of using the flipped classroom model is the varying methods of implementation. Since there is no consistent model for a flipped classroom, the concept rarely means the same thing to instructors in terms of implementation [5]. O’Flaherty and Phillips also share this concern; they concluded from their cross-disciplinary scoping review that “there is a danger that educators renewing their curriculum may not fully understand the pedagogy of how to effectively translate the flipped class into practice…there appears to be some misunderstanding of the key elements necessary for successful flipping” [11]. According to them, these key elements include the consideration of ways to improve students’ engagement both inside and outside of class and the facilitation of critical thinking [11]. Another aspect of successful implementation frequently discussed in the literature is student-professor communication. Robert Talbert believes this is what differentiates successful and unsuccessful flipped classrooms. According to him, professors should not only clearly state why the classroom is inverted to students both early in the semester and throughout the course but also listen to students’ opinions about the course structure and be accommodating to needed changes in the course to improve their learning [12]. Many articles emphasize that professors who decide to try the flipped classroom technique should be mindful of their method of implementation and be willing to change their methodology during the semester if necessary.
\section{AVAILABLE RESEARCH}
    \hspace{4mm} Many of the research studies note that there is not yet enough quantitative research on the flipped classroom technique compared with the traditional lecture method to make definite conclusions on its efficacy. Love et al. mention that though the flipped classroom technique has shown potential for improving student learning and interest in STEM areas, most research and discussion has been “more anecdotal than data driven” [2]. However, the amount of total research on the method has greatly increased in the past few years. In Flipped Learning: A Guide for Higher Education Faculty, Robert Talbert includes a graph that shows the number of peer-reviewed published articles that contain the words “flipped classroom,” “inverted classroom,” or “flipped learning” in their titles or abstracts published each year between 2011 and 2015[6]. His data demonstrates an approximate increase of 250 percent per year in the quantity of peer-reviewed articles published on flipped learning and “suggests that flipped learning is not just a fad or a buzzword but a pedagogical model being seriously explored and examined by a rapidly growing segment of the academic population” [6]. He also notes that “as professionals in higher education, we often want more than just good explanation and compelling stories when we think about a new idea. We want evidence, preferably from sound scholarship and reputable research. After all, we are talking about a pedagogical model that seems experimental and novel.” According to Talbert, professors who are considering implementing the flipped model and those who have already switched want to learn more about the reasons why the method yields positive results: “[if] a professor might be switching to flipped learning from a traditional approach to teaching, she might like to have some reasonable expectation of why this way of running a class should work as well or better than the methods higher education has adopted for a very long time. On the other hand, a professor who has used flipped learning before can benefit from knowing why certain things about it are potentially transformative for student learning” [6].

{\setlength{\parskip}{3mm}
       In addition to the growing number of articles published on the flipped classroom method within specific subjects at certain colleges, there have been studies that consider the effectiveness of the flipped classroom model across many disciplines in higher education. In 2015, Jaqueline O’Flaherty and Craig Phillips published a scoping review of the use of the flipped classroom in higher education. In their research, they initially identified 1084 articles involving the flipped classroom and narrowed their study to 28 articles by using inclusion criteria such as population and sample, study focus, and article type [11]. There have also been several large scale article reviews in STEM disciplines other than math; in 2017, Karabulut-Ilgu et al. published a systematic review of research on flipped learning within engineering education [17]. This study reviewed 62 articles on the flipped classroom technique that were published between 2000 and 2015 [17]. In addition, in 2016 Huber and Werner published a review of the literature on flipped STEM classrooms. This study analyzed 58 peer-reviewed research studies that focused on the flipped classroom technique within STEM disciplines in higher education and identified both positive and negative themes in the articles regarding students’ perception, engagement, self-efficacy achievement, and development of graduate attributes [18]. Further, Rahman et al. published a meta-analysis of the influences of the flipped classroom in 2014 that analyzed the results of 15 articles, most of which focused on STEM disciplines within higher education [19].
       {\setlength{\parskip}{3mm}
Although there are systemic reviews of the research on the flipped classroom technique within higher education, STEM disciplines, and engineering education, there is not yet a systemic review of the flipped classroom technique within higher math education in particular. Robert Talbert mentions that such a study would be an important next step in research on the flipped classroom technique in order to measure normalized gains in student learning from pre-test and post-test scores between flipped and traditional sections of specific math courses [20].
    {\setlength{\parskip}{3mm}
\section{QUANTITATIVE RESEARCH}
      \hspace{4mm} Although there are not yet large-scale or meta-analysis studies on the effectiveness of the flipped classroom method for mathematics, there are many studies providing quantitative research that yield an indication of the effectiveness of the flipped classroom method within certain math classes and populations of students. These studies compared student exam results from flipped and traditional class sections. The quantitative data from these articles is summarized in this chapter and key aspects of the findings of these studies are introduced. Undergraduate math courses that have been studied include Topics in Calculus, Calculus I, Calculus II, Calculus III, Transition-to-Proof, and Linear Algebra.
\subsection{Topics in Calculus}
      \hspace{4mm} In a small-scale study, Bill Blubaugh compared the effectiveness of the flipped classroom technique to a semi-traditional method in two Topics in Calculus courses he taught in Fall 2014. The students in both sections were primarily business majors, and there were 37 students in the flipped classroom and 32 students in the semi-traditional classroom. The format of the semi-traditional class included a short lecture and group work. The most significant quantitative findings of this research were the results of the first test: the class average for the flipped class was 88.76 percent while the average for the semi-traditional class was 48.88 percent. He notes, however, that the reason for the large difference between averages is debatable and not necessarily only due to the two differing instructional methods. On average, students in the flipped class were older and had achieved higher grades previously, which could have helped their comprehension and exam performance levels relative to the students in the semi-traditional class [21]. The averages for both the second test and final for the flipped and semi-traditional classes were much closer to one another. Blubaugh’s exam research is summarized below.
     \\ \\
\includegraphics[scale=.7]{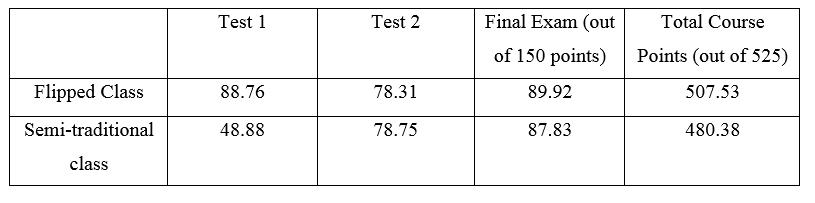}
\\
{\setlength{\parskip}{3mm} \subsection{Calculus 1}
     \hspace{4mm} In a relatively larger study at Sungkyunkwan University in Korea, Karjanto and Simon studied the effectiveness of the flipped classroom technique in a single variable calculus course over three semesters, from Spring 2015 to Spring 2016. To gather quantitative research, they compared the average of midterm and final exam scores as well as the distribution of the letter grades given to students at the end of the semester. They compared four different types of classroom settings: a fully flipped class using instructor-made videos, a fully flipped class using videos from Khan Academy, a class that was only flipped for a single topic (Introduction to Differential Questions), and a traditional lecture class. They concluded there was no statistically significant difference in terms of letter grades between the four types of classroom settings. However, there was a statistically significant difference between the fully flipped section and the partially flipped section in terms of assessment averages. They note that there were significantly more students retaking the course in the single topic flipped class than the other classes, which likely contributed to this class having higher assessment averages and average letter grades [9]. There were 59 students in the flipped class that used the instructor-made videos, 85 students in the flipped class that used the Khan Academy videos, 81 students in the single topic flipped class, and 85 students in the traditional class. The quantitative data of Karjanto and Simon’s study is summarized below [9]. Note that assessment averages are the average of the midterm exam and the final exam, and in terms of the value of the letter grade in calculating students’ GPAs, Sungkyunkwan University uses the following:
A+ = 4.5 A = 4.0 B+ = 3.5 B = 3.0 C+ = 2.5 C = 2.0 D+ = 1.5 D = 1.0 and F = 0. \\ \\
\includegraphics[scale=.7]{blubaugh.png}\\
\subsection{Calculus 2}

Jean McGivney-Burelle and Fei Xue studied the implementation of the flipped classroom technique in Calculus II at the University of Hartford in Spring 2012. During this semester, the same instructor taught two sections of Calculus II and used the traditional lecture method during the entire semester for one of the sections and switched to the flipped classroom method in the other section for a unit that is typically challenging to students [16]. There were 29 students in the section that used the traditional lecture format for the entire semester and 31 students in the section that switched to the flipped class format. McGivney-Burelle and Xue found that students’ grades on exams and homework were higher in the flipped classroom section than the traditional section for this unit. The exam and homework grade data from their study is summarized below, and they concluded that the flipped class method has “potential to work well across a range of mathematics” and note that they “look forward to developing more flipped units” for their calculus courses [16].\\
\includegraphics[scale=.7]{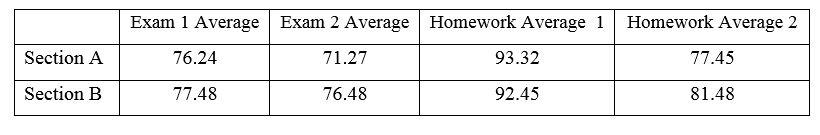}\\
Note that Section A used the traditional lecture format for the entire semester and Section B switched to the flipped class format for the material that was assessed in Exam 2 and Homework 2. Therefore, the averages for Exam 1 and Homework 1 can be considered a baseline for Section A compared to Section B.
\subsection{Calculus 3}
     Wasserman et al. studied the effectiveness of the flipped classroom technique compared with the traditional lecture style in regards to both students’ performance and perceptions in Calculus III courses [5]. This study was completed at Southern Methodist University over two semesters, Fall 2012 and Spring 2013, in order to gather data from a larger number of students. A total of 74 students in the flipped classes and a total of 77 students in the traditional classes participated in the study. There are a few important points to note in considering the data from this study. First, prior to Test 1, both sections used the traditional (non-flipped) method. After Test 1, one class changed to the flipped method. Two professors taught the sections involved in the study, each teaching one section of the Calculus III course both semesters. The professors used identical lecture notes, homework assignments, and assessments. The flipped class switched to using video lectures for the more procedural components and the more conceptual components continued to be taught in class. This study also gathered and analyzed much demographic data between student populations in the traditional and flipped classes, including each student’s major, Calculus II grade, SAT math score, year in college, age, and gender. It concluded that the two populations were similar, with no significant differences in demographics between the traditional and flipped classes overall [5]. The students in both sections also performed similarly on the first test. Wasserman et al. concluded that similar scores on Test 1 along with the similar demographic data provide a baseline for comparing the effectiveness of both instructional methods in the results of the second and third exams. The findings of their research showed that students in the flipped class performed slightly to moderately better than students in the traditional class [5]. Furthermore, Wasserman et al. analyzed the average of the conceptual questions as well as the procedural questions and concluded that the higher averages for students in the flipped classroom were due to their better performance of the conceptual test questions than the students in the traditional class. They found that students in the flipped class did not perform significantly better or worse on the procedural questions than students in the traditional classroom. Combined from both semesters, the overall test averages along with the averages from the conceptual and procedural questions from this study are summarized below.\\
\includegraphics[scale=.7]{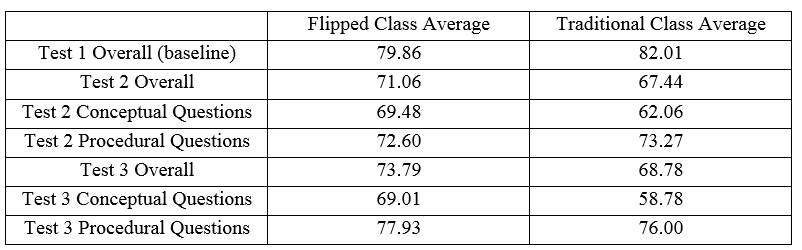}\\
     This study also notes that in comparing the data gathered in Fall 2012 and Spring 2013, there are larger differences between the exam averages for the flipped and traditional classes for the Spring 2013 semester (with students in the flipped class performing better). Wasserman et al. attribute this to improvement in implementation of the flipped method, due to the prior experience with using the technique. They interpret this as evidence for potential improvements of flipped classroom implementation over time [5]. Specifically, they think these improvements were likely due to more productive class time, their creation of better and more concise video lectures, and the incorporation of clicker quizzes to hold students more accountable for outside-of-class work. Overall, this study provides many ideas for further research, particularly in its methodologies and efforts to isolate the flipped classroom technique as a single variable in assessing student learning.
\subsection{Transition to Proof}
     \hspace{4mm} Robert Talbert studied the benefits of implementing the flipped class technique in a transition-to-proof course for math majors at Grand Valley State University. Talbert mentions that a goal of transition-to-proof courses is for students to be able to focus on specifics in subsequent proof-based math courses and not on the actual process of proof writing. However, this is not always the learning outcome of students who complete these classes, and many students struggle with basic proof-related tasks in later courses. Talbert references a study of 61 students in a transition-to-proof course that found that none of them could consistently rewrite an informal mathematical statement as a formal and logically equivalent statement [22]. Further, another study that tracked students’ and professors’ eye movements as they read various proofs found that undergraduates focused significantly less on the actual arguments and significantly more on the superficial aspects than the professional mathematicians [23].
     Talbert notes that there are inherently many pedagogical challenges concomitant with courses that serve as a transition from computational to proof-based math. In this course, students transition from a primarily computational view of mathematics to a primarily conceptual view. This transition requires students to acquire self-regulated learning strategies that are not always developed and practiced in lower-level math classes [12]. In addition, transition-to-proof courses require students to internalize sociomathematical norms. According to Yackel and Cobb, sociomathematical norms are “normative aspects of mathematical discussions that are specific to students’ mathematical activity” [24]. Prior to taking a transition-to-proof course, students typically have little experience with professional mathematics and the accompanying appropriate professional norms. Talbert mentions that “most college and high school students do not know what a proof is or what it is supposed to achieve” [12].\\ \\
    {\setlength{\parskip}{3mm}  Due to these pedagogical challenges present in transition-to-proof courses, Talbert studied the benefits of implementing the flipped classroom method in two sections of the course. Talbert compared students’ grade distributions from two flipped class sections with a total of 39 students to students’ grade distributions from traditional classes from Fall 2007 to Fall 2012 with a total of 939 students. Talbert concludes that implementing the flipped classroom technique in this class did not significantly change course grades. However, in the flipped section, there was a marginal increase in top grades and a marginal decrease in non-passing grades. Talbert notes that as instructors using the flipped classroom technique make improvements in their implementation and instruction, there could be larger decreases in students who withdraw from the course or receive non-passing grades [12]. This idea of increased improvement in future semesters supports the arguments of Wasserman et al. The data from Talbert’s transition-to-proof study is summarized below [12]. Note that “W” refers to students who withdrew from the course.
\\ \\
\includegraphics[scale=.7]{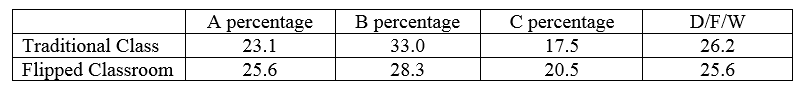}\\

\subsection{Linear Algebra}
     Betty Love et al. studied student learning in a flipped classroom section compared with a traditional lecture section of an applied linear algebra course during Spring 2012 [2]. Their study consisted of 27 students in the flipped classroom section and 28 students in the traditional lecture section. Students in both sections were primarily math, computer science, and engineering majors. The study measured student learning by comparing students’ exam results from three midterm exams and a comprehensive final. They note that “of particular interest is the analysis of how the students progressed between exams, since the benefit of the flipped class approach may not be evident early in the course but should increase as the course progresses” [2]. They found that the average change in score between the first and second midterms was significantly greater in the flipped classroom section than in the traditional lecture section. The average change in score between the second and third midterm was also significantly greater in the flipped classroom section than in the traditional lecture section. Students in the flipped classroom section also performed slightly better on the final exam; the average raw score for the students in the flipped section was an 89.5, and the average in the traditional lecture section was 87.4 [2].
\section{QUANTITATIVE RESEARCH}
Much of the literature published on the flipped classroom method also includes qualitative research primarily centered around three interconnected ideas: student perceptions, student self-efficacy, and increased peer interaction. The results of individual studies with key findings are summarized in this chapter. Many of these studies yield ideas for further research on the effectiveness of the flipped classroom model.
\subsection{Student Perceptions}

     Sparking students’ interest in math is particularly important in introductory courses. Students who have the ability and interest to continue in math are more likely to persist if they have a positive experience in their initial college math courses. This is important not only for students on an individual basis but also on a national scale with the increased need for students in math and math-based fields in the United States. Thus, it is important to consider students’ perceptions when implementing a new teaching method, and a significant number of researchers have focused on students’ perceptions to the flipped classroom implementation or included it in their studies. \\ \\
     Wasserman et al. emphasize the importance of examining students’ perceptions of pedagogical methods. They mention that research establishing the interaction between students’ performance and students’ perceptions is increasingly needed due to the widespread use of the flipped classroom technique [5]. In their study, Wasserman et al. analyzed student perceptions in the flipped and traditional classroom sections of Calculus III through anonymous surveys given to students in both semesters that the study was completed. On average, students in the flipped class sections were more likely to agree that they communicated during class regularly than students in the traditional sections. However, students in the traditional class were more likely to agree that class time was used effectively [5]. The student perception data from both semesters of the Wasserman et al. study, as component mean comparisons, is summarized below. \\
\includegraphics[scale=.7]{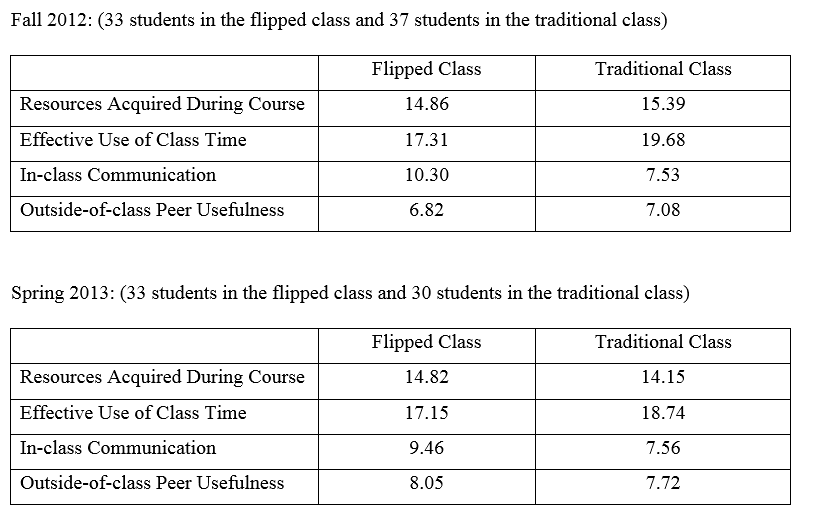}\\
     Further, students in this study who participated in a focus group mentioned that the high difficulty level of the problems that they were assigned in class caused class time to be unproductive, as they had to wait for the instructor’s explanation if their group didn’t know how to solve it. Additionally, some students mentioned that they thought the flipped classroom method caused them to miss the instructor’s expertise, as they were primarily working with other students who were novices at the concepts [5]. \\ \\
     Love et al. also studied students’ perceptions in flipped classrooms through student surveys and found positive results. In their study, they found that 74 percent  of students in the flipped classroom had a positive attitude about the flipped class approach, and students in the flipped class section found the video lectures “significantly more helpful” than the students in the lecture section found the lectures. In fact, 96 percent  of students in the flipped section thought the videos helped them to learn the material. The results from the student surveys showed that approximately one-third of the students who watched the lectures watched them more than once, and one student commented that “Many times I have wished to pause or rewind a live lecture. You can do that with a video without disrupting the flow of the class. Also, watching a second time really reinforces the concepts” [2]. In addition, over 74 percent of students in the flipped section of this study agreed that working problems on the board helped them to better remember course concepts and was more enjoyable than in a traditional lecture. The students’ comments showed similarly positive perceptions: “The more interactive environment of the course held my attention and helped me stay focused” and “This style of teaching kept my attention…it is more fun to work with others on the board than it is working by yourself on paper” [2]. Love et al. also concluded from surveys that students in the flipped class section completed the course with “a greater perception that linear algebra is relevant to their career” than students in the traditional section. \\ \\
     Robert Talbert also gathered qualitative data from students in his flipped classroom section of a transition-to-proof course. Though only sixteen students answered open-ended questions about experiences with the flipped classroom, their responses provide some insight into student perceptions of the flipped class method. In a question asking students to describe their overall experiences with the course, thirteen students had positive comments about their experience with the course format, and only two of the comments were negative. Talbert included one: “I felt like I was thrown off the boat and expected to float…Book was a little helpful, but was very frustrated when asking for help and I couldn’t have questions answered” [12]. He mentions that this comment illustrates that students in flipped classes need to have much support for the course, in multiple ways and accessible in different places, and this support should be “advertised to students on a continuous basis.” Specifically, he notes that including a brief explanation on each class handout that states the way in which students can request help while in different contexts of their learning (while they are watching the lecture video, studying outside of class, etc.) “significantly stemmed” comments of this nature. \\ \\
     Talbert also asked students through the survey whether they believed that their experiences in the class improved their ability to learn new content on their own. Twelve students agreed, three disagreed, and one was unsure. Some students agreeing to this statement included comments that offer insight into their perceptions to the course structure. One mentioned that “The inverted classroom setting has taught me how to do things more on my own as opposed to just listening to what the professor says,” and another mentioned that “We worked through the lecture on our own and it gave us the opportunity to work through our problems first, without instantly being rescued. It was frustrating at times but I guess overall I have benefited from it” [12]. When asked how their attitudes and strategies about learning math changed over the course of the semester, eleven students reported a positive change and two reported no change due to having already engaged in flipped classroom methodology before taking the course. Two of the students who self-reported a positive change in their attitudes and learning strategies elaborated: “I would like to take another class in this format, I think it was helpful…I think it is a much more viable way of learning topics like the ones in the class” and “I am looking forward to taking math classes of similar styles in the future” [12]. \\ \\
     In the final question of the survey, Talbert asked students how their attitudes and strategies about learning non-mathematical subjects changed over the semester. Four out of twelve students mentioned a positive change in their learning skills, specifically in improved time-management and discipline in their other classes. The remaining eight students reported no change. Talbert concluded that these numbers suggest “a modest indication that students built global self-regulated learning skills, but it more strongly indicated that a single inverted course may not catalyze a global change in students’ learning strategies all by itself. Instructors who wish to use an inverted classroom as a means of improving general study strategies and self-regulated learning behaviors in students will probably need to make a point of addressing the notion of self-regulated learning in class and lead students to think about how their problem-solving skills generalize to other areas of study.” Talbert’s qualitative findings from these student surveys are summarized below [12]. Note that not all sixteen students answered every survey question.\\ \\
\includegraphics[scale=.7]{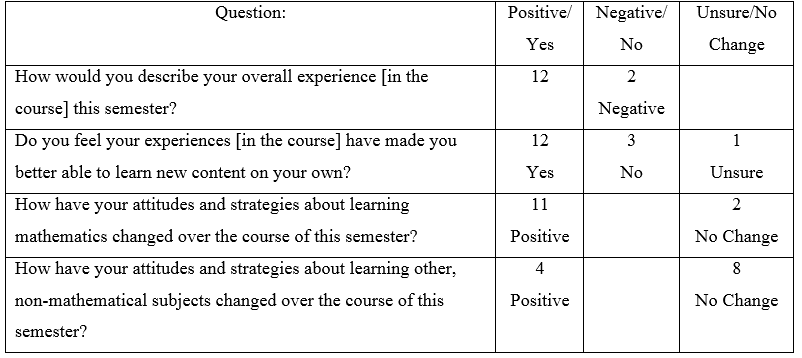}\\
\\
     Talbert notes that many students will “experience the inverted classroom as a major culture shock and will rebel against it” [12]. He elaborates on this idea in another article, noting that many students’ initial experience with the flipped classroom is “troubling and uncomfortable” because of their lack of familiarity with it and their belief that their professor’s job is to lecture [20]. However, he believes that “student discomfort over the lack of in-class lecturing can give way to meaningful discussions about the nature of higher education and real progress toward guid[ing] students to becoming self-regulated lifelong learners” [20].\\
    \\ In McGivney-Burelle and Xue’s study comparing flipped classroom and traditional sections of a Calculus II course, they conducted focus group interviews with four students to assess their perceptions of the flipped classroom method. All four students had strong preferences for flipped classroom instruction over traditional instruction. The aspects of the flipped classroom method that the students most preferred were the ability to watch videos outside of class and the way class time was spent [16]. One student mentioned that during class, she felt more comfortable asking questions “because there is less of a time constraint. I would not ask questions when he was just lecturing because…I don’t want to hold up the whole class.” Another student mentioned that watching the videos and doing examples is “more of a self-tester” and “boosts your confidence.” Another student called the videos a “stress-free visual method of learning” [16]. To improve flipped classroom implementation, the students recommended that instructors spend a few minutes at the beginning of class answering any questions about the video lectures since some students wished they were able to ask their instructors questions while watching the videos [16].\\
\subsection{Student Self-Efficacy}

Some researchers conclude that the flipped classroom technique can improve students’ self-efficacy by making them more independent learners. Albert Bandura defines self-efficacy as “people’s beliefs about their capabilities to produce designated levels of performance that exercise influence over events that affect their lives” [26]. Further, Bandura mentions that “a strong sense of efficacy enhances human accomplishment and personal well-being in many ways. People with high assurance in their capabilities approach difficult tasks as challenges to be mastered rather than as threats to be avoided. Such an efficacious outlook fosters interest and deep engrossment in activities” [26]. Lori Ogden studied student perceptions in flipped College Algebra courses over the course of three semesters. She concluded that at the end of the semester, students responded that they were more confident in their ability to learn math. Even students who mentioned not liking math felt they were “more capable and no longer fear mathematics” [25]. Hall and Ponton mention that research has shown that two key elements for success in mathematics are perceived ability and previous mathematical performance [27]. Ogden mentions that because of this established link between students’ self-efficacy and achievement in math courses, it is “critical that instructors adopt pedagogical techniques that mediate the effects of low self-efficacy” [25]. The results of Ogden’s study suggest that flipped classroom implementation may help students to feel more in control of their own learning [25].\\ \\
Similarly, Robert Talbert credits the flipped classroom technique for helping students to become more self-regulated learners. He thinks that “lecturing may inadvertently create a dependence of the student upon the lecturer, thereby preventing students from developing fully as self-regulated learners” [12]. According to Talbert, many students “have developed a nearly invulnerable sense of learned helplessness through their primary and secondary education that conflates learning with lecture and with high numerical scores on tests.” Talbert mentions that these ideas can contribute to student success overall since “helping students emerge as competent, confident, self-regulating learners is hard work but is the primary job of higher education. The inverted classroom shows promise of catalyzing progress toward this goal” [12].\\ \\
Matthew Voight studied students’ attitudes on the flipped classroom technique in four pre-calculus II sections that were taught at a large research university in the Midwest. Two of the sections used the flipped classroom method, and two of the sections used the traditional lecture method. His study included 427 student responses to pre-course and post-course surveys, and he found significant differences between perceived student experiences in the flipped sections and the traditional sections. Students in the flipped sections self-reported significantly higher ratings for cognitive, affective, and social learning gains as well as math confidence and collaborative strategies for problem solving [28]. These results indicate that the implementation of the flipped classroom could improve students’ self-efficacy. \\ \\
     Another idea related to student self-efficacy is student attendance in flipped versus tradition sections of classes. O’Flaherty and Phillips reference one study within their meta-analysis that found that attendance increased by 30 percent to 80 percent when the flipped classroom model was used rather than the traditional lecture model [11]. On the other hand, Blubaugh mentions that students in the flipped class section missed an average 6.43 classes and students in the semi-traditional section missed an average 4.74 per student [21], though the small size of this study should be noted. There is not yet enough data regarding student attendance in flipped classrooms versus traditional classrooms to make definite conclusions, so it would be an interesting topic to include in further research.
     Another student survey question in Blubaugh’s study is related to student self-efficacy: “How much time do you spend each week preparing for class (studying, reading, homework, office hours, tutor lab, videos, etc)?” He found that, on average, students in the flipped class section spent a little more time studying than students in the traditional section. However, the small sample size should again be noted, as a total of only 48 students answered this question. The results are summarized below.\\ \\
\includegraphics[scale=.7]{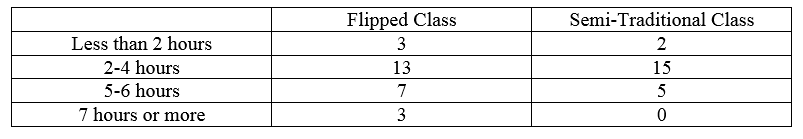}\\
\\     Countless articles emphasize the importance of improving retention within STEM majors, and research has begun to link aspects of flipped classroom implementation, specifically improved student self-efficacy, to improved retention within STEM majors. In 2014, Boise State University implemented a reform of instruction in all of their calculus courses due to weaknesses in the calculus sequence, including an average pass rate of 51 percent . The reformed version of Calculus 1 at Boise State University included aspects of the flipped classroom technique; the majority of class time was allocated to students working on assignments in small groups designed to target specific learning goals [29]. Bullock et al. studied how this reform of calculus courses influenced retention within STEM majors. They analyzed retention data from 1175 students who took calculus classes under the original pedagogical methods and 1177 students who took calculus courses with the reformed teaching methods. They found that the new version of calculus increased STEM retention by 3.3 percent overall. It had a more significant impact on two groups: STEM retention for women increased by 9.1 percent and STEM retention for underrepresented minorities increased by 9.4 percent . Bullock et al. notes that students’ self-efficacy is a critical element strongly associated with retention in STEM disciplines [29]. They reference a 2016 study that found that women are 1.5 times more likely to leave STEM disciplines after calculus compared to men, and this study attributes this to lack of mathematical confidence [30]. Ellis et al. found that when comparing women and men with above-average mathematical abilities and preparedness, women begin and end semesters with significantly lower confidence in mathematics than men. From their research, Ellis et al. concludes “if women persisted in STEM at the same rate as men starting in Calculus I, the number of women entering the STEM workforce would increase by 75 percent” [30]. Bullock et al. mention that their study suggests potential for improvements in students’ self-efficacy with the implementation of the reformed version of calculus since the retention gap between men and women closed virtually entirely [29]. \\
\subsection{Peer Interaction}
     Helping to build connections between students of similar majors is beneficial to students not only socially but also academically, as students in more interactive classes are more likely to study together than students in less interactive ones. Research has shown that implementation of the flipped classroom technique increases peer interaction and in turn helps students to build connections and improve their collaborative learning abilities. \\ \\
     From a survey given to students at the end of the semester, Love et al. found that while 56 percent of students in the flipped section agreed that they were more comfortable talking to classmates in that particular class than other math courses they had taken, only 21 percent of students in the traditional section of the class agreed [2]. Love et al. mentions that the increased amount of peer-interaction that the flipped classroom technique encompasses can contribute to the development of social networks for students. They note that this is particularly helpful for the retention of women and minority students within math [2].\\ \\
     Bullock et. al also emphasize the importance of working to build a community learning atmosphere within classes. They refer to this idea as a need for belongingness and mention that a need to belong and to be accepted within a community of peers drives peoples’ behaviors and choices [29]. Further, this need increases in adverse or stressful environments such as high stakes academic ones. In his article titled “Ingroup Experts and Peers as Social Vaccines Who Inoculate the Self-Concept: The Stereotype Inoculation Model,” Nilanjana Dasgupta mentions that “the experience of being in a numeric minority in high-stakes academic environments where stereotypes are in the air may reduce individuals’ self-efficacy or confidence in their own abilities, especially in the face of difficulty, even if their actual performance is objectively the same as majority-group members” [31]. Because of this, the increased peer interaction within flipped classes could improve students’ sense of belongingness, especially if they are members of underrepresented groups. \\ \\
     If this theory holds that the flipped classroom model increases retention of students within math, it could be a significant advantage of the flipped classroom model due to the increased need for students to pursue STEM majors in the US. In 2012, the President’s Council of Advisors on Science and Technology, an advisory group composed of leading scientists and engineers in the U.S., called for an increase of 33 percent in the number of bachelor’s degrees in science, technology, engineering, and math completed per year. It found that fewer than 40 percent of students in the United States who begin college with an interest in a STEM fields graduate with a STEM degree [32]. The disparity is ever greater for students in underrepresented groups; only 20 percent of underrepresented minority students who were initially interested in STEM field finish with a STEM degree [7]. Further, one of the five overarching recommendations of the report was “to catalyze widespread adoption of empirically validated teaching practices,” which links to the flipped classroom model. The report also emphasizes the importance of creating “an atmosphere of a community of STEM learners” [32]. This attests to the potential power of the increased peer-interaction that the flipped classroom model provides. \\ \\
\section{CONCLUSIONS}
     The flipped classroom technique has shown great potential for enhancing students’ learning and experiences within their classes, including in college mathematics. Research on the subject is growing exponentially and shows that students in flipped classroom sections of various undergraduate math courses typically perform slightly to moderately better than students in traditional lecture courses. With the expansion of implementation and refinement of implementation methods, there is potential for more substantial gains. The benefits in terms of students’ experiences are even more apparent, with research indicating improved student perception, improved student self-efficacy, and increased peer interaction within flipped classes. Positive student experiences within undergraduate math courses are vital for retention within math and math-related areas of study.\\ \\
     With increased interest in and application of the flipped classroom method, instructors should be cognizant of their methods of implementation. Used carefully, technology can augment student learning in flipped classes, and video lectures in some form are an important component that students prefer. Student-instructor communication in flipped classrooms is key. Students should be informed of their instructor’s reasoning for structuring the course in the flipped class method, so that students are more accepting of the method. Instructors should clearly communicate their expectations for pre-class assignments and the benefits of completing them. Additionally, instructors should be willing to adapt aspects of course structure if students find them ineffective or have ideas for improved implementation. Instructors can incorporate aspects of the flipped classroom into their instruction without completely flipping their classes. If possible, math departments could give students the option whether to enroll in flipped or non-flipped sections of classes in order to better match student learning styles. \\ \\
     In reviewing the literature on the flipped classroom technique, a clear theme is a continued need for additional research as the technique is applied in a greater variety of math courses and implemented in different methods. Used carefully, the flipped classroom technique could greatly improve student learning and interest in undergraduate math classes as well as contribute to the self-regulated learning goals of higher education.\\  \\
\section{ACKNOWLEDGEMENTS}
Special thanks to Dr. Furuzan Ozbek, project advisor, for all of her guidance, help in finding resources, and explanations her own experiences with applying the flipped classroom technique. Thank you to Dr. Betty Love and Dr. Neal Grandgenett for sending additional information from their study.

{\centering
{\setlength{\parskip}{3mm}\section*{References}}}
1. Bishop, Jacob and Matthew Verleger. 2013. “The Flipped Classroom: A Survey of the Research.” American Society for Engineering Education Conference and Exposition. Atlanta, GA. \\ \\
2. Love, Betty, Angie Hodge, Neal Grandgennet, and Andrew W. Swift. 2014. “Student Learning and Perceptions in a Flipped Linear Algebra Course” International Journal of Mathematical Education in Science and Technology 45(3): 317-324. \\ \\
3. Felder, Richard. 1993. “Reaching the Second Tier: Learning and Teaching Styles in College Science Education.” J. College Science Teaching 23(5): 286-290. \\ \\
4. Lord, Susan, and Michelle Camacho. 2007. “Effective Teaching Practices: Preliminary Analysis of Engineering Educators.” 37th ASEE/IEEE Frontiers in Education Conference. Milwaukee, WI. \\ \\
5. Wasserman, Nicholas H., Christa Quint, Scott A. Norris, and Thomas Car. 2015. “Exploring Flipped Classroom Instruction in Calculus III.” International Journal of Science and Math Education 15:545-568. \\ \\
6. Talbert, Robert. Flipped Learning: A Guide for Higher Education Faculty. 2017. Sterling, VA: Stylus. \\ \\
7. Freeman, Scott, Sarah L. Eddy, Miles McDonough, Michelle K. Smith, Nnadozie Okoroafor, Hannah Jordt, and Mary P. Wenderoth. 2014. “Active Learning Increases Students Performance in Science, Engineering, and Mathematics.” Proceedings of the National Academy of Sciences 111(23): 8410-8415. \\ \\
8. Brame, Cynthia. 2013. “Flipping the classroom.” Vanderbilt University Center for Teaching. \\ \\
9. Karjanto, Natanael. and L. Simon. 2017. “Lecture at Home and Homework at School: Flipped Class 101 in an English-Medium Instruction Single Variable Calculus Course.” \\ \\
10. “Definition of Flipped Learning.” 2014. Flipped Learning Network.  \\ \\
11. O’Flaherty, Jacqueline and Craig Phillips. 2015. “The Use of Flipped Classrooms in Higher Education: A Scoping Review.” Internet and Higher Education 25: 85-95. \\ \\
12. Talbert, Robert. 2015. “Inverting the Transition-to-Proof Classroom.” Problems, Resources, and Issues in Mathematics Undergraduate Studies 25(8): 614-626. \\ \\
13. McGrath, Dominic, Anthea Groessler, Esther Fink, Carl Reidsema, and Lydia Kavanagh. 2017. “Technology in the Flipped Classroom.” The Flipped Classroom: Practice and Practices in Higher Education. Singapore: Springer. pp. 37-56. \\ \\
14. Williams, Cassandra and John Siegfried. 2016. “Investigating Student-Learning Gains from Video Lessons in a Flipped Calculus Course.” Special Interest Group of the Mathematical Association of America on Research in Undergraduate Mathematics Education. \\ \\
15. Armstrong, Patricia. “Bloom’s Taxonomy.” Vanderbilt University Center for Teaching. \\ \\
16. McGivney-Burelle, Jean and Fei Xue. 2013. “Flipping Calculus.” Problems, Resources, and Issues in Mathematics Undergraduate Studies 23(5): 477-486. \\ \\
17. Karabulut-Ilgu, Aliye, Nadia Jaramillo Cherrez, and Charles T. Jahren. 2017. “A Systemic Review of Research on the Flipped Learning Method in Engineering Education.” British Journal of Educational Technology. \\ \\
18. Huber, Elaine and Ashleigh Werner. 2016. “A Review of the Literature on Flipping the STEM Classroom.” Australian Society for Computers in Learning in Tertiary Education Conference. Adelaide, Australia. \\ \\
19. Rahman, Azlina, Baharuddin Aris, Hasnah Mohamed, and Norasykin Mohd Zaid. 2014. “The Influences of the Flipped Classroom: A Meta Analysis.” IEEE 6th International Conference on Engineering Education. Kuala Lumpur, Malaysia. \\ \\
20. Talbert, Robert. 2014. “Inverting the Linear Algebra Classroom.” Problems, Resources, and Issues in Mathematics Undergraduate Studies 24(5): 361-374. \\ \\
21. Blubaugh, Bill. 2015. “A Comparison between a Flipped-Learning Calculus Class and a Semi-Traditional Class.” 27th International Conference on Technology in Collegiate Mathematics. Las Vegas, NV. \\ \\
22. 22. Selden, John and Annie Selden. 1995. “Unpacking the Logic of Mathematical Statements.” Educational Studies in Mathematics 29(2): 123-151. \\ \\
23. Inglis, Matthew and Lara Alcock. 2012. “Expert and Novice Approaches to Reading Mathematical Proofs.” Journal for Research in Mathematics Education 43(4): 358-390. \\ \\
24. Yackel, Erna and Paul Cobb. 1996. “Sociomathematical Norms, Argumentation, and Autonomy in Mathematics.” Journal for Research in Mathematics Education 27(4): 458-477. \\ \\
25. Ogden, Lori. 2015. “Student Perceptions of the Flipped Classroom Technique in College Algebra.”  Problems, Resources, and Issues in Mathematics Undergraduate Studies 25(9-10): 782-791.\\ \\
26. Bandura, Albert. 1994. “Self-Efficacy.” Encyclopedia of Human Behavior. 4: 71-81. \\ \\
27. Hall, Michael and Michael K. Ponton. 2005. “Mathematics Self-Efficacy of College Freshman.” Journal of Developmental Education 28(3): 26-28. \\ \\
28. Voight, Matthew. 2016. “Examining Student Attitudes and Mathematical Knowledge Inside the Flipped Classroom Experience.” Conference on Research in Undergraduate Mathematics Education. Pittsburgh, PA. \\ \\
29. Bullock, Doug, Janet Callahan, and Jocelyn Cullers. 2017. “Calculus Reform: Increasing STEM Retention and Post-Requisite Course Success While Closing the Retention Gap for Women and Underrepresented Minority Students.” American Society for Engineering Education Annual Conference. Columbus, OH. \\ \\
30. Ellis, Jessica, Bailey Fosdick, and Chris Rasmussen. 2016. “Women 1.5 Times More Likely to Leave STEM Pipeline after Calculus Compared to Men: Lack of Mathematical Confidence a Potential Culprit.” Public Library of Science One 11(7). \\ \\
31. Dasgupta, Nilanjana. 2011. “Ingroup Experts and Peers as Social Vaccines Who Inoculate the Self-Concept: The Stereotype Inoculation Model.” An International Journal for the Advancement of Psychological Theory. 22(4): 231-246. \\ \\
32. Gates, S.J., Jo Handelsman, G.P. Lapage, and Chad Mirken. 2012. “PCAST STEM Undergraduate Working Group.” Engage to Excel: Producing One Million Additional College Graduates with Degrees in Science, Technology, Engineering, and Mathematics. Office of the President, Washington.

\end{document}